\documentclass[letterpaper,11pt]{article}

 \title{A subcopula based dependence measure}
 \author{Arturo Erdely}
 \date{\small{Facultad de Estudios Superiores Acatl\'an \\
              Universidad Nacional Aut\'onoma de M\'exico \\
							\texttt{arturo.erdely@comunidad.unam.mx}\\}}
 
 \usepackage{amsthm}
 \usepackage{amssymb}
 \usepackage{amsmath}
 \usepackage{enumerate}
 \usepackage{graphicx}

\begin{document}
 
\maketitle

\begin{abstract}
\noindent A dependence measure for arbitrary type pairs of random variables is proposed and analyzed, which in the particular case where both random variables are continuous turns out to be a concordance measure. Also, a sample version of the proposed dependence measure based on the empirical subcopula is provided, along with an R package to perform the corresponding calculations.
\end{abstract}


\section{Introduction}

If $(X,Y)$ is a bivariate random vector with joint probability distribution $F_{X,Y}(x,y)=P(X\leq x, Y\leq y),$ the outstanding theorem by Sklar (1959) ensures that there exists a unique functional relationship $S$ between $F_{X,Y}$ and its marginal univariate probability distribution functions $F_X(x)=P(X\leq x)$ and $F_Y(y)=P(Y\leq y),$ such that:
\begin{equation}\label{Sklar}
F_{X,Y}(x,y)\,=\,S(F_X(x),F_Y(y))\,,\quad x,y\in\overline{\mathbb{R}}=[-\infty,+\infty].
\end{equation}
Since the ranges of $F_{X,Y}, F_X,$ and $F_Y$ are subsets of the unit interval $\mathbb{I}=[0,1]$ which at least include $0$ and $1$ then $S$ is a function with domain $\text{Ran }F_X\times\text{Ran }F_Y\subseteq\mathbb{I}^{\,2}$ and range a subset of $\mathbb{I}$ which at least includes $0$ and $1.$ As an immediate consequence of (\ref{Sklar}) we obtain:
\begin{itemize}
  \item[a)] $0=F_{X,Y}(-\infty,y)=S(F_X(-\infty),F_Y(y))=S(0,v)$ where $v\in\text{Ran }F_Y,$ and analogously $S(u,0)=0$ for $u\in\text{Ran }F_X.$
	\item[b)] $F_Y(y)=F_{X,Y}(+\infty,y)=S(F_X(+\infty),F_Y(y))=S(1,v)$ where $v=F_Y(y),$ and analogously $S(u,1)=u$ where $u=F_X(x)$ for some $x\in\overline{\mathbb{R}}.$
	\item[c)] $0\leq P(x_1<X\leq x_2\,,\,y_1<Y\leq y_2)=F_{X,Y}(x_2,y_2)-F_{X,Y}(x_2,y_1)-F_{X,Y}(x_1,y_2)+F_{X,Y}(x_1,y_1)$ and therefore by (\ref{Sklar}) we have that $S(u_2,v_2)-S(u_2,v_1)-S(u_1,v_2)+S(u_1,v_1)\geq 0$ where $u_i=F_X(x_i)$ and $v_i=F_Y(y_i)$ for $i=1,2.$
\end{itemize}

\medskip

\noindent\textbf{Definition 1} A bivariate \textit{subcopula} (or 2-subcopula) is a function $S:D_1\times D_2\rightarrow\mathbb{I},$ where $\{0,1\}\subseteq D_i\subseteq\mathbb{I}$ $(i=1,2),$ such that for all $u,v\in\mathbb{I}:$
\begin{itemize}
  \item[a)] $S(u,0)=0=S(0,v)\hspace{0.3mm};$
	\item[b)] $S(u,1)=u\,$ and $\,S(1,v)=v\hspace{0.3mm};$
	\item[c)] $S(u_2,v_2)-S(u_2,v_1)-S(u_1,v_2)+S(u_1,v_1)\geq 0\,$ where $\,u_1\leq u_2\,$ and $\,v_1\leq v_2\hspace{0.3mm}.$
\end{itemize}

\medskip

\noindent Therefore the unique functional relationship $S$ in (\ref{Sklar}) is a subcopula. In the particular case when the domain of a bivariate subcopula is $D_1\times D_2=\mathbb{I}^{\,2}$ then it is called bivariate \textit{copula} (or 2-copula). This will be the case when both $X$ and $Y$ are continuous random variables, but in any other case $D_1\times D_2$ will be a proper subset of $\mathbb{I}^{\,2}.$ Any subcopula which is not a copula may be extended to a copula in a non-unique way, see for example Lemma 2.3.5 in Nelsen (2006). Every subcopula $S$ is bounded by the Fr\'echet-Hoeffding bounds:
\begin{equation}\label{FH}
W(u,v)\leq S(u,v)\leq M(u,v)
\end{equation}
where $W(u,v)=\max\{u+v-1,0\}$ and $M(u,v)=\min\{u,v\}$ are copulas which may be restricted to subcopulas with the same domain as subcopula $S,$ denoted by $W_S$ and $M_S,$ respectively. Recalling that $X$ and $Y$ are independent random variables (of any kind) if and only if $F_{X,Y}(x,y)=F_X(x)F_Y(y),$ the unique underlying subcopula for such random vector $(X,Y)$ according to (\ref{Sklar}) would be $S(u,v)=uv$ where the domain of $S$ would be $\text{Ran }F_X\times\text{Ran }F_Y.$ It is common to use the notation $\Pi(u,v)=uv$ which is a copula that may also be restricted to any subcopula domain, for example $\Pi_S$ as in the notation introduced before. As an immediate consequence of Theorems 2.5.4 and 2.5.5 in Nelsen (2006) we obtain the following:

\bigskip

\noindent\textbf{Corollary 1} \textit{Let $(X,Y)$ be a random vector such that $Y=g(X)$ for some function $g,$ and let $S$ be its unique underlying subcopula according to (\ref{Sklar}).}
  \begin{itemize}
	  \item[a)] \textit{$S=M_S$ if and only if $g$ is almost surely nondecreasing on $\text{Ran}\,X.$}
  	\item[b)] \textit{$S=W_S$ if and only if $g$ is almost surely nonincreasing on $\text{Ran}\,X.$}
	\end{itemize}
	
\smallskip

\noindent In the particular case $X$ and $Y$ are continuous random variables, as explained in Nelsen (2006):
\begin{quote}
  \textsl{When $X$ and $Y$ are continuous, the support of their joint distribution function can have no horizontal or vertical segments, and in this case it is common to say that ``$Y$ is almost surely [a strictly] increasing function of $X$'' if and only if the copula of $X$ and $Y$ is $M;$ and ``$Y$ is almost surely a [strictly] decreasing function of $X$'' if and only if the copula of $X$ and $Y$ is $W.$}
\end{quote}

\noindent As discussed in Ne\v{s}lehov\'a (2007) and Genest and Ne\v{s}lehov\'a (2007) for continuous random variables many dependence concepts and measures of association can be expressed in terms of the unique underlying copula only and thus independently from the marginal distributions. This interrelationship fails as soon as there are discontinuities in the marginal distribution functions: the possibility of ties that results from atoms in the probability distributions invalidates various familiar relations that lie at the root of copula theory in the continuous case, and so neither the axiomatic definition for a concordance measure by Scarsini (1984) nor the use of the concordance function is clear. Moreover, as stated in Ne\v{s}lehov\'a (2007):
\begin{quote}
  \textsl{The fact that marginal distributions functions take influence upon the dependence structure is characteristic for non-continuous distributions. In the case of concordance measures, this ``nuissance'' causes difficulties: the measures typically do not reach the bounds $\pm 1$ for countermonotonic and comonotonic marginals.}
\end{quote}

\smallskip

\noindent The probabilistic definitions of popular concordance measures (such as Kendall or Spearman) do not account for ties, so modified versions of their theoretical and empirical definitions are needed, see: Denuit and Lambert (2005), Ne\v{s}lehov\'a (2007), and Genest \textit{et al.}(2014); but the way to define them is non-unique since they are based on non-unique extensions of subcopulas to copulas. The main contribution of the present work is to propose a dependence measure based directly on the unique underlying subcopula, regardless of the random variable types in a bivariate random vector $(X,Y).$ Since a copula $C$ is a particular case of subcopula, such proposal turns out to be a concordance measure for a pair of continuous random variables that is related to the $L_{\infty}$ distance between $C$ and $\Pi$ in an similar way as Spearman's concordance measure is related to the $L_1$ distance between $C$ and $\Pi$ known as Schweizer and Wolff (1981) dependence measure.

\section{A monotone dependence measure}

Let $S:D_1\times D_2\rightarrow\mathbb{I}$ be the unique underlying subcopula for a random vector $(X,Y)$ of arbitrary type random variables accordingly to (\ref{Sklar}), where $\{0,1\}^2\subseteq D_1\times D_2=\text{Ran }F_X\times\text{Ran }F_Y\subseteq\mathbb{I}^{\,2}.$

\bigskip

\noindent\textbf{Definition 2} (Adapted from Lehmann 1966 and Nelsen 2006). Two random variables $X$ and $Y$ will be called \textit{positively quadrant dependent} (PQD) if $P(X\leq x\,,\,Y\leq y)\geq P(X\leq x)P(Y\leq y),$ which by (\ref{Sklar}) is equivalent to $S(u,v)\geq\Pi_S(u,v)$ for all $(u,v)\in D_1\times D_2.$ Negative quadrant dependence (NQD) is defined analogously by reversing the sense of the inequalities, that is $S(u,v)\leq\Pi_S(u,v)$ for all $(u,v)\in D_1\times D_2.$ 

\bigskip

\noindent\textbf{Proposition 1} \textit{Let $\mathcal{A}$ be the set of all bivariate subcopulas. The function $d:\mathcal{A}\rightarrow\mathbb{R}$ defined by}
\begin{equation}\label{dS}
   d(S)\,:=\,\sup_{\text{Dom}\,S}\{S-\Pi_S\}\,-\,\sup_{\text{Dom}\,S}\{\Pi_S-S\}\,,\qquad S\in\mathcal{A}
\end{equation}
\textit{has the following properties:}
\begin{itemize}
  \item[a)] $d(\Pi_S)=0\,;$
	\item[b)] $-\frac{1}{4}\leq d(W_S)\leq 0\leq d(M_S)\leq\frac{1}{4}\,;$
	\item[c)] $d(W_S)\leq d(S)\leq d(M_S)\,;$
	\item[d)] $|d(S)|=\frac{1}{4}$ \textit{if and only if} $(\frac{1}{2},\frac{1}{2})\in\text{Dom}\,S\,;$
	\item[e)] \textit{if $X$ and $Y$ are PQD (respectively NQD) then $d(S)\geq 0$ (respectively $d(S)\leq 0$)}\,;
	\item[f)] \textit{if $S_1,S_2\in\mathcal{A}$ such that $\text{Dom}\,S_1=D=\text{Dom}\,S_2$ and $S_1\preceq S_2$ then $d(S_1)\leq d(S_2).$}
\end{itemize}

\noindent\textbf{Proof:}
\begin{itemize}
  \item[a)] It follows immediately from the definition. 
	\item[b)] Elementary calculations can show that for a function $h:\mathbb{I}^{\,2}\rightarrow\mathbb{R}$ defined by $h(u,v):=M(u,v)-uv\geq 0$ we have that $\max h=\frac{1}{4}=h(\frac{1}{2},\frac{1}{2})$ and $h(u,v)<\frac{1}{4}$ for $(u,v)\neq(\frac{1}{2},\frac{1}{2})\,,$ and therefore
	 $$0\,\stackrel{\text{by (\ref{FH})}}{\leq}\,d(M_S) \,=\, \sup_{\text{Dom}\,S}\{M_S-\Pi_S\} \,-\,0\quad\left\{ \begin{array}{cc}
		                                                     = \frac{1}{4} & \text{ if } (\frac{1}{2},\frac{1}{2})\in\text{Dom}\,S \\
																												       { }     &     { } \\
														    												 < \frac{1}{4} & \text{ if } (\frac{1}{2},\frac{1}{2})\notin\text{Dom}\,S
									      															 \end{array}\right.$$
and an analogous result follows using $g(u,v):=W(u,v)-uv\leq 0$ since $\min g=-\frac{1}{4}=g(\frac{1}{2},\frac{1}{2})$ and $g(u,v)>-\frac{1}{4}$ for $(u,v)\neq(\frac{1}{2},\frac{1}{2})\,.$ 
  \item[c)] From (\ref{FH}) we get $S-\Pi_S\preceq M_S-\Pi_S$ and then
  $$d(S)\,\leq\,\sup_{\text{Dom}\,S}\{S-\Pi_S\}\,\leq\,\sup_{\text{Dom}\,S}\{M_S-\Pi_S\}-0\,=\,\sup_{\text{Dom}\,S}\{M_S-\Pi_S\}-\sup_{\text{Dom}\,S}\{\Pi_S-M_S\}\,=\,d(M_S)$$
and an analogous reasoning leads to $d(W_S)\leq d(S).$ 
  \item[d)] It follows immediately from the arguments to prove b) and c). 
	\item[e)] If $X$ and $Y$ are PQD then their underlying unique subcopula $S\succeq\Pi_S$ and therefore $\sup_{\text{Dom}\,S}\{\Pi_S-S\}=0$ so we conclude from the definition that $d(S)\geq 0.$ Case NQD is analogous.
	\item[f)] If $S_1\leq S_2$ then $\sup_D\{S_2-\Pi_D\}\geq\sup_D\{S_1-\Pi_D\}$ and $-\sup_D\{\Pi_D-S_2\}\geq-\sup_D\{\Pi_D-S_1\}\,,$ and by adding left and right sides of these two inequalities we obtain $d(S_2)\geq d(S_1).\qquad_{\pmb{\square}}$
\end{itemize}

\bigskip

\noindent\textbf{Definition 3} The \textit{monotone dependence measure} for arbitrary type random variables $X$ and $Y$ with underlying subcopula $S$ will be denoted and defined as:
  $$\mu_{X,Y}\,\equiv\mu(S)\,:=\,\left\{\begin{array}{cc}
	                                          d(S)/d(M_S) & \text{ if } d(S)\geq 0 \text{  and } M_S\neq\Pi_S\,, \\
																			  		-d(S)/d(W_S) & \text{ if } d(S)\leq 0 \text{  and } W_S\neq\Pi_S\,, \\
																				  	      0      & \text{ if } W_S=\Pi_S=M_S\,.
	                                      \end{array}\right.$$
Notice that $W_S=\Pi_S=M_S$ will occur if and only if Dom$\,S=\{0,1\}\times\{0,1\},$ and that would be the case of a pair of constant random variables, so in what follows we will not consider this trivial case. The expression \textit{monotone dependence measure} is not new, it is being used in a similar way as in Cifarelli \textit{et al.}(1996).

\bigskip

\noindent\textbf{Theorem 1} \textit{The monotone dependence measure satisfies the following properties:}
\begin{itemize}
  \item[a)] $\mu_{X,Y}$ \textit{is defined for every pair of arbitrary type random variables;}
	\item[b)] $-1\leq\mu_{X,Y}\leq +1\,,$ $\mu_{X,X}=+1\,,$ \textit{and} $\mu_{X,-X}=-1\,;$
	\item[c)] $\mu_{X,Y}=\mu_{Y,X}\,;$
	\item[d)] \textit{if $X,Y$ are independent then $\mu_{X,Y}=0\,;$}
	\item[e)] if $\mu_{X,Y}\neq 0$ then $\mu_{X,Y}$ has opposite sign to $\mu_{X,-Y}$ and $\mu_{-X,Y}\,;$
	\item[f)] \textit{if $S_1,S_2$ are subcopulas such that $\text{Dom}\,S_1=D=\text{Dom}\,S_2$ and $S_1\preceq S_2$ then $\mu(S_1)\leq\mu(S_2)\,;$}
	\item[g)] \textit{if $P[Y=\varphi(X)]=1$ with $\varphi$ nondecreasing (respectively nonincreasing) then $\mu_{X,Y}=+1$ (respectively $\mu_{X,Y}=-1)\,;$}
	\item[h)] \textit{if $X$ and $Y$ are PQD (respectively NQD) then $\mu_{X,Y}\geq 0$ (respectively $\mu_{X,Y}\leq 0$).}
\end{itemize}

\noindent\textbf{Proof:}
\begin{itemize}
  \item[a)] An immediate consequence of Sklar's theorem (\ref{Sklar}) since $\mu_{X,Y}$ is defined in terms of the unique underlying subcopula $S.$
	\item[b)] If $d(S)\geq 0$ then by Definition 3 we have that $\mu_{X,Y}=d(S)/d(M_S)$ and applying Proposition 1 c) we get $\mu_{X,Y}\leq +1.$ Similarly, if $d(S)\leq 0$ then $\mu_{X,Y}=-d(S)/d(W_S)$ and by Proposition 1 b) and c) we now get $\mu_{X,Y}\geq -1.$ Now by Corollary 1 using $Y=X$ we have that the underlying subcopula for $(X,X)$ is $S=M_S$ so $d(S)=d(M_S)\geq 0$ and therefore $\mu_{X,X}=d(M_S)/d(M_S)=+1.$ Similarly, using $Y=-X$ we have that the underlying subcopula for $(X,-X)$ is $S=W_S$ so $d(S)=d(W_S)\leq 0$ and therefore $\mu_{X,-X}=-d(W_S)/d(W_S)=-1.$
	\item[c)] Straightforward by applying Sklar's theorem to the fact that $F_{X,Y}(x,y)=P(\{X\leq x\}\cap\{Y\leq y\})=P(\{Y\leq y\}\cap\{X\leq x\})=F_{Y,X}(y,x).$
	\item[d)] If $X$ and $Y$ are independent then $F_{X,Y}(x,y)=F_X(x)F_Y(y)$ and by Sklar's theorem their unique underlying subcopula is $S=\Pi_S$ so by Proposition 1a we have $d(S)=0$ and therefore $\mu_{X,Y}=0.$
	\item[e)] It will suffice to prove that $d(S_{X,-Y})=-d(S_{X,Y})$ where $S_{X,Y}$ and $S_{X,-Y}$ are the unique underlying copulas for $(X,Y)$ and $(X,-Y),$ respectively.
	    \begin{equation}\label{eq1}
	       F_{-Y}(y)\,=\,P(-Y\leq y)\,=\,P(Y\geq -y)\,=\,1-F_Y(-y)+P(Y=-y)\,=\,1-F_Y((-y)^{-})
			\end{equation} where $F_Y((-y)^{-})=\lim_{z\rightarrow(-y)^{-}}F_Y(z)$ is a left-hand limit at $-y,$ and so
			\begin{equation}\label{eq2}
			   \text{Ran}\,F_{-Y}\,=\,\{F_{-Y}(y):y\in\overline{\mathbb{R}}\}\,=\,\{1-F_Y(y^{-}):y\in\overline{\mathbb{R}}\}\,.
			\end{equation}
			Let $S_{X,Y}:\text{Ran}\,F_X\times\text{Ran}\,F_Y\rightarrow\mathbb{I}$ be the unique underlying subcopula for $(X,Y),$ and also let $S_{X,-Y}:\text{Ran}\,F_X\times\text{Ran}\,F_{-Y}\rightarrow\mathbb{I}$ be the unique underlying subcopula for $(X,-Y).$ Then:
			\begin{eqnarray*}
			  S_{X,-Y}\big(F_X(x),F_{-Y}(y)\big)\,\stackrel{\text{Sklar}}{=}\,F_{X,-Y}(x,y) &=& P(X\leq x,-Y\leq y) = P(X\leq x, Y\geq -y) \\
				                                     &=& P(X\leq x) - \lim_{z\rightarrow(-y)^{-}}P(X\leq x, Y\leq z) \\
																						 &=& F_X(x) - F_{X,Y}\big(x,(-y)^{-}\big) \\
																						 &\stackrel{\text{Sklar}}{=}& F_X(x) - S_{X,Y}\big(F_X(x),F_Y((-y)^{-})\big) \\
																						 &\stackrel{\text{(\ref{eq1})}}{=}& F_X(x) - S_{X,Y}\big(F_X(x),1-F_{-Y}(y)\big).
			\end{eqnarray*}
			If we define $u:=F_X(x)$ and $v:=F_{-Y}(y)$ in this last result we get:
			\begin{equation}\label{eq3}
			  S_{X,-Y}(u,v)\,=\,u\,-\,S_{X,Y}(u,1-v)
			\end{equation}
			where $u\in\text{Ran}\,F_X$ and $v$ must satisfy:
			\begin{equation}\label{eq4}
			  v\in\text{Ran}\,F_{-Y}\quad\text{ and }\quad 1-v\in\text{Ran}\,F_Y\,,
			\end{equation}
			which implies that $v\in\text{Ran}\,F_{-Y}=\{1-F_Y(y-):y\in\overline{\mathbb{R}}\}$ and that $v\in\text{Ran}(1-F_Y):=\{1-F_Y(y):y\in\overline{\mathbb{R}}\},$ where $\text{Ran}\,F_{-Y}=\text{Ran}(1-F_Y)=\mathbb{I}\,$ if $Y$ is a continuous random variable, otherwise the symmetric difference $\text{Ran}\,F_{-Y}\bigtriangleup\text{Ran}(1-F_Y)$ is at most countable, in which case both sides of (\ref{eq3}) can be properly defined over the domain $\text{Ran}\,X\times D$ where the set $D$ is the closure of $\text{Ran}\,F_{-Y}\cap\text{Ran}(1-F_Y),$ by taking adequate limits. Finally:
			\begin{eqnarray*}
			  d(S_{X,-Y}) &\,=\,& \sup_{(u,v)\in D}\{S_{X,-Y}(u,v)-uv\}\,-\,\sup_{(u,v)\in D}\{uv-S_{X,-Y}(u,v)\} \\
			      &\,\stackrel{\text{(\ref{eq3})}}{=}\,& \sup_{(u,v)\in D}\{u-S_{X,Y}(u,1-v)-uv\}\,-\,\sup_{(u,v)\in D}\{uv-u+S_{X,Y}(u,1-v)\} \\
						&\,=\,& \sup_{(u,v)\in D}\{u(1-v)-S_{X,Y}(u,1-v)\}\,-\,\sup_{(u,v)\in D}\{S_{X,Y}(u,1-v)-u(1-v)\} \\
						&\,=\,& -d(S_{X,Y}).
			\end{eqnarray*}
	\item[f)] An immediate consequence of Proposition 1f and Definition 3.
	\item[g)] Applying Corollary 1, if $\varphi$ is almost surely nondecreasing then the underlying subcopula for $(X,Y)$ is $S=M_S$ so $d(S)=d(M_S)\geq 0$ and therefore $\mu_{X,Y}=d(M_S)/d(M_S)=+1.$ Similarly, if $\varphi$ is almost surely nonincreasing then the underlying subcopula for $(X,Y)$ is $S=W_S$ so $d(S)=d(W_S)\leq 0$ and therefore $\mu_{X,Y}=-d(W_S)/d(W_S)=-1.$	
  \item[h)] An  immediate consequence of Proposition 1e and Definition 3. $\qquad_{\pmb{\square}}$
\end{itemize}

\bigskip

\noindent\textbf{Corollary 2} \textit{If $X$ and $Y$ are continuous random variables with unique underlying copula $C$ then:}
\begin{itemize}
  \item[a)] $\mu_{X,Y}\,\equiv\,\mu(C)\,=\,4\big(\max_{\,\mathbb{I}^{\,2}}\{C-\Pi\}\,-\,\max_{\,\mathbb{I}^{\,2}}\{\Pi-C\}\big)\,;$
	\item[b)] $\mu$ \textit{is a measure of concordance.}
\end{itemize}
\noindent\textbf{Proof:}
\begin{itemize}
  \item[a)] An immediate consequence from Proposition 1 and Definition 3 since $d(M)=\frac{1}{4}=d(W),$ and the fact that $C$ is a continuous function with domain the compact set $\mathbb{I}^{\,2}.$
	\item[b)] Theorem 1 includes all the properties for a measure of concordance required by Definition 5.1.7 in Nelsen (2006) except the following one: if $\{(X_n,Y_n)\}$ is a sequence of continuous random variables  with copulas $C_n,$ and if $\{C_n\}$ converges pointwise to $C,$ then $\lim_{\,n\rightarrow\infty}\mu(C_n)=\mu(C),$ but this is straightforward to prove since all copulas are continuous with domain the compact set $\mathbb{I}^{\,2},$ and so we may exchange maximum and limit. $\qquad_{\pmb{\square}}$
\end{itemize}

\noindent From Nelsen (2006), in the particular case of copulas the $L_{\infty}$ distance between $C$ and $\Pi$ is given by
\begin{equation}\label{supdist}
  \Lambda(C)\,=\,4\sup_{\mathbb{I}^2}\big|C\,-\,\Pi\big|,
\end{equation}
and therefore if $C$ is a member of a totally ordered (with respect to the concordance ordering $\preceq$) copula family that includes $\Pi$ then $\Lambda(C)=|\mu(C)|,$ in a similar way as $\sigma(C)=|\rho(C)|$ where $\sigma$ is Schweizer and Wolff (1981) dependence measure and $\rho$ is Spearman's concordance measure:
\begin{equation}\label{SpearmanSchweizer}
  \sigma(C) = 12\!\int\!\!\!\int_{\mathbb{I}^{\,2}}\!|C(u,v)-uv|\,dudv\,,\qquad\rho(C) = 12\!\int\!\!\!\int_{\mathbb{I}^{\,2}}\![C(u,v)-uv]\,dudv\,.
\end{equation}

\section{Examples}

Let $(X,Y)$ be a random vector of arbitrary type random variables with joint distribution function
\begin{equation}\label{Sklar2}
  F_{XY}(x,y\,|\,\alpha,\theta_1,\theta_2)\,=\,S_{\alpha}\big(\,F_X(x\,|\,\theta_1)\,,\,F_Y(y\,|\,\theta_2)\,\big)
\end{equation}
where $(\alpha,\theta_1,\theta_2)$ belongs to some parametric space, and where accordingly to (\ref{Sklar}) we have that the function $S_{\alpha}:\text{Ran}\,F_X\times\text{Ran}\,F_Y\rightarrow\mathbb{I}$ is the unique underlying subcopula with parameter $\alpha\,;$ and $\theta_1$ and $\theta_2$ are marginal parameters of $X$ and $Y,$ respectively. If both $X$ and $Y$ are continuous random variables then $\text{Ran}\,F_X=\mathbb{I}=\text{Ran}\,F_Y,$ which implies that the domain of $S_{\alpha}$ is $\mathbb{I}^{\,2}$ and therefore $S_{\alpha}$ would be, in fact, a copula. In this particular case, the value of $\mu_{X,Y}$ will be only a function of the subcopula parameter $\alpha.$ In case one of the random variables is non-continuous, say $X,$ then $\text{Ran}\,F_X$ is a proper subset of $\mathbb{I}$ which depends on the marginal parameter $\theta_1$ and therefore $\mu_{X,Y}$ will be a function of $\alpha$ and (possibly) $\theta_1.$ And of course if both random variables are non-continuous the value of $\mu_{X,Y}$ will be a function of $\alpha$ and (possibly) of $\theta_1$ and/or $\theta_2.$

\bigskip\medskip

\noindent\textbf{Example 1} Consider a bivariate random vector $(X,Y)$ where $X$ and $Y$ are Bernoulli random variables with parameters $0<\theta_1<1$ and $0<\theta_2<1,$ respectively, and dependence parameter $\alpha = P(X=1,Y=1).$ We may summarize its joint probability mass function $P(X=x,Y=y)$ as:
\begin{equation}\label{pmfBer2D}
  \begin{array}{|c|c|c|c|} \hline
	  P(X=x,Y=y) &             Y=0             &       Y=1       & P(X=x)     \\ \hline
		      X=0  & 1+\alpha -\theta_1-\theta_2 & \theta_2-\alpha & 1-\theta_1 \\
				  X=1  &       \theta_1-\alpha       &      \alpha     & \theta_1   \\ \hline
        P(Y=y) &         1-\theta_2				   &      \theta_2   &    { }     \\ \hline
  \end{array}
\end{equation} 
where by Fr\'echet-Hoeffding bounds $\max\{\theta_1+\theta_2-1,0\}\leq\alpha\leq\min\{\theta_1,\theta_2\},$ $X$ and $Y$ are independent if and only if $\alpha=\theta_1\theta_2,$ and PQD/NQD if $\alpha\geq\theta_1\theta_2$ or $\alpha\leq\theta_1\theta_2,$ respectively. In this case the unique underlying subcopula $S:D_1\times D_2\rightarrow\mathbb{I}$ is determined by $D_1=\{0,1-\theta_1,1\},$ $D_2=\{0,1-\theta_2,1\},$ and $S(1-\theta_1,1-\theta_2)=1+\alpha-\theta_1-\theta_2$ since the other 8 subcopula values are determined by boundary conditions a) and b) in Definition 1:
\begin{equation}\label{subcopulaEx1}
  \begin{array}{|l|c|c|c|} \hline
	            S(u,v) & v = 0 & v = 1-\theta_2                 & v = 1          \\ \hline
		           u = 0 & 0 & 0                          & 0          \\ \hline
		  u = 1-\theta_1 & 0 & 1+\alpha-\theta_1-\theta_2 & 1-\theta_1 \\ \hline
					     u = 1 & 0 & 1-\theta_2                 & 1          \\ \hline
	\end{array}
\end{equation}
Then $d(S)=\alpha-\theta_1\theta_2$ and
\begin{equation}\label{dMS}
  d(M_S)\,=\,\left\{ \begin{array}{cc}
	                         \theta_2(1-\theta_1)\,, & \theta_1\geq\theta_2 \\
				     							           { }           &          { }         \\
													 \theta_1(1-\theta_2)\,, & \theta_1\leq\theta_2 \\
											\end{array}\right.
\end{equation}
\begin{equation}\label{dWS}
  -d(W_S)\,=\,\left\{ \begin{array}{cc}
	                                 \theta_1\theta_2\,, & \theta_2\leq 1-\theta_1 \\
														           	  { }          &           { }           \\
									         (1-\theta_1)(1-\theta_2)\,, & \theta_2\geq 1-\theta_1 \\
											\end{array}\right.
\end{equation}
Therefore applying Definition 3:
\begin{equation}\label{muXYejem1}
  \mu_{X,Y}\,=\,\begin{cases}
	                  (\alpha-\theta_1\theta_2)/\theta_2(1-\theta_1)\,, & \theta_2\leq\theta_1\,,\,\alpha\geq\theta_1\theta_2 \\
										(\alpha-\theta_1\theta_2)/\theta_1(1-\theta_2)\,, & \theta_2\geq\theta_1\,,\,\alpha\geq\theta_1\theta_2 \\
										(\alpha-\theta_1\theta_2)/\theta_1\theta_2\,, & \theta_2\leq 1-\theta_1\,,\,\alpha\leq\theta_1\theta_2 \\
										(\alpha-\theta_1\theta_2)/(1-\theta_1)(1-\theta_2)\,, & \theta_2\geq 1-\theta_1\,,\,\alpha\leq\theta_1\theta_2
	              \end{cases}
\end{equation}
It is straightforward to obtain $\text{Cov}(X,Y)=\alpha-\theta_1\theta_2=d(S),$ and therefore Pearson's correlation coefficient $r_{X,Y}=(\alpha-\theta_1\theta_2)/\sqrt{\theta_1(1-\theta_1)\theta_2(1-\theta_2)}.$ In Table 1 it is compared $\mu_{X,Y}$ versus $r_{X,Y}$ under extreme values of $\alpha$ and independence:
\begin{table}[h]
  \begin{center}
    \begin{tabular}{|c|c|c|c|c|} \hline
                 $\alpha$                 &      condition         &       interpretation    & $\mu_{X,Y}$ & $r_{X,Y}$ \\ \hline
	        $\min\{\theta_1,\theta_2\}$     & $\theta_1\neq\theta_2$ &    $P(Y\geq X) = 1$     &    $+1$     &  $< +1$   \\
					$\min\{\theta_1,\theta_2\}$     &  $\theta_1=\theta_2$   &      $P(Y = X) = 1$     &    $+1$     &   $+1$    \\
					    $\theta_1\theta_2$          &         none           & $X$ and $Y$ independent &     $0$     &    $0$    \\
					$\max\{\theta_1+\theta_2-1,0\}$ & $\theta_2>1-\theta_1$  &   $P(Y\geq 1-X) = 1$    &    $-1$     &  $> -1$   \\
					$\max\{\theta_1+\theta_2-1,0\}$ & $\theta_2<1-\theta_1$  &   $P(Y\leq 1-X) = 1$    &    $-1$     &  $> -1$   \\
					$\max\{\theta_1+\theta_2-1,0\}$ & $\theta_2=1-\theta_1$  &   $P(Y = 1-X) = 1$      &    $-1$     &   $-1$    \\ \hline
    \end{tabular}
    \caption{Comparing monotone dependence measure versus Pearson's correlation in Example 1.}
    \label{ejemplo1}
  \end{center}
\end{table}

\bigskip

\noindent\textbf{Example 2} Let $(X,Y)$ be a random vector with joint probability distribution function:
\begin{equation}\label{conjunta2}
  F_{X,Y}(x,y\,|\,\alpha,\theta)\,=\,\alpha\min\{1-x^{-1},1-(1-\theta)^{\lfloor y \rfloor}\}+(1-\alpha)(1-x^{-1})[1-(1-\theta)^{\lfloor y \rfloor}]\,,\quad x>1,\,y\geq 1
\end{equation}
with parameters $0\leq\alpha\leq 1$ and $0<\theta<1,$ and where $\lfloor y \rfloor$ stands for the maximum integer less than or equal to $y.$ By marginalization it is straighforward to verify that $X$ is a continuous random variable Pareto$(1,1)$ and $Y$ is a discrete Geometric$(\theta)$ random variable since $F_X(x)=F_{X,Y}(x,+\infty)=1-x^{-1},$ $x>1,$ and $F_Y(y\,|\,\theta)=F_{X,Y}(+\infty, y)=1-(1-\theta)^{\lfloor y \rfloor},$ $y\geq 1,$ and by (\ref{Sklar2}) it is obtained:
\begin{equation}\label{SklarEjem2}
  F_{X,Y}(x,y\,|\,\alpha,\theta)\,=\,S_{\alpha}\big(\,F_X(x)\,,\,F_Y(y\,|\,\theta)\,\big)\,,\quad x>1,\,y\geq 1
\end{equation}
with underlying subcopula $S_{\alpha}:\text{Ran}\,F_X\times\text{Ran}\,F_Y\rightarrow\mathbb{I}^{\,2},$ where $\text{Ran}\,F_X=\mathbb{I}$ and $\text{Ran}\,F_Y=\{1-(1-\theta)^k:k=0,1,\ldots\}\cup\{1\}\subset\mathbb{I},$ given by:
\begin{equation}\label{subcopEjem2}
  S_{\alpha}(u,v)\,=\,\alpha\min\{u,v\}\,+\,(1-\alpha)uv\,=\,\alpha M(u,v)\,+\,(1-\alpha)\Pi(u,v)\,,
\end{equation}
that is, subcopula $S_\alpha$ is a convex linear combination of copulas $M$ and $\Pi$ restricted to $\text{Ran}\,F_X\times\text{Ran}\,F_Y,$ where $S_0=\Pi$ and $S_1=M.$ Recalling (\ref{FH}) we have that $\Pi\leq M$ which in this example implies $S_\alpha\geq\Pi$ and therefore:
\begin{eqnarray*}
  d(S_{\alpha}) &\,=\,& \sup_{\text{Dom}\,S_{\alpha}}\{\,S_{\alpha}(u,v)\,-\,\Pi_{S_{\alpha}}(u,v)\,\}\,-\,0 \\
	              &\,=\,& \sup_{\text{Dom}\,S_{\alpha}}\{\,\alpha\hspace{0.3mm}[\,M(u,v)-\Pi(u,v)\,]\,\}\,=\,\alpha\sup_{\text{Dom}\,S_{\alpha}}\{\,M(u,v)-\Pi(u,v)\,\} \\
								&\,=\,& \alpha\hspace{0.3mm}d(M_{S_{\alpha}})
\end{eqnarray*}
and consequently $\mu_{X,Y}=d(S_{\alpha})/d(M_{S_{\alpha}})=\alpha.$ It should be noticed that despite $Y$ is a discrete random variable its parameter $\theta$ does not have an influence on $\mu_{X,Y}$ in this case. In this example it is not possible to calculate Pearson's correlation since first moment of $X$ does not exist.$\qquad_{\pmb{\square}}$

\bigskip\medskip

\begin{center}
\begin{figure}[h]
\includegraphics[width=17cm, keepaspectratio]{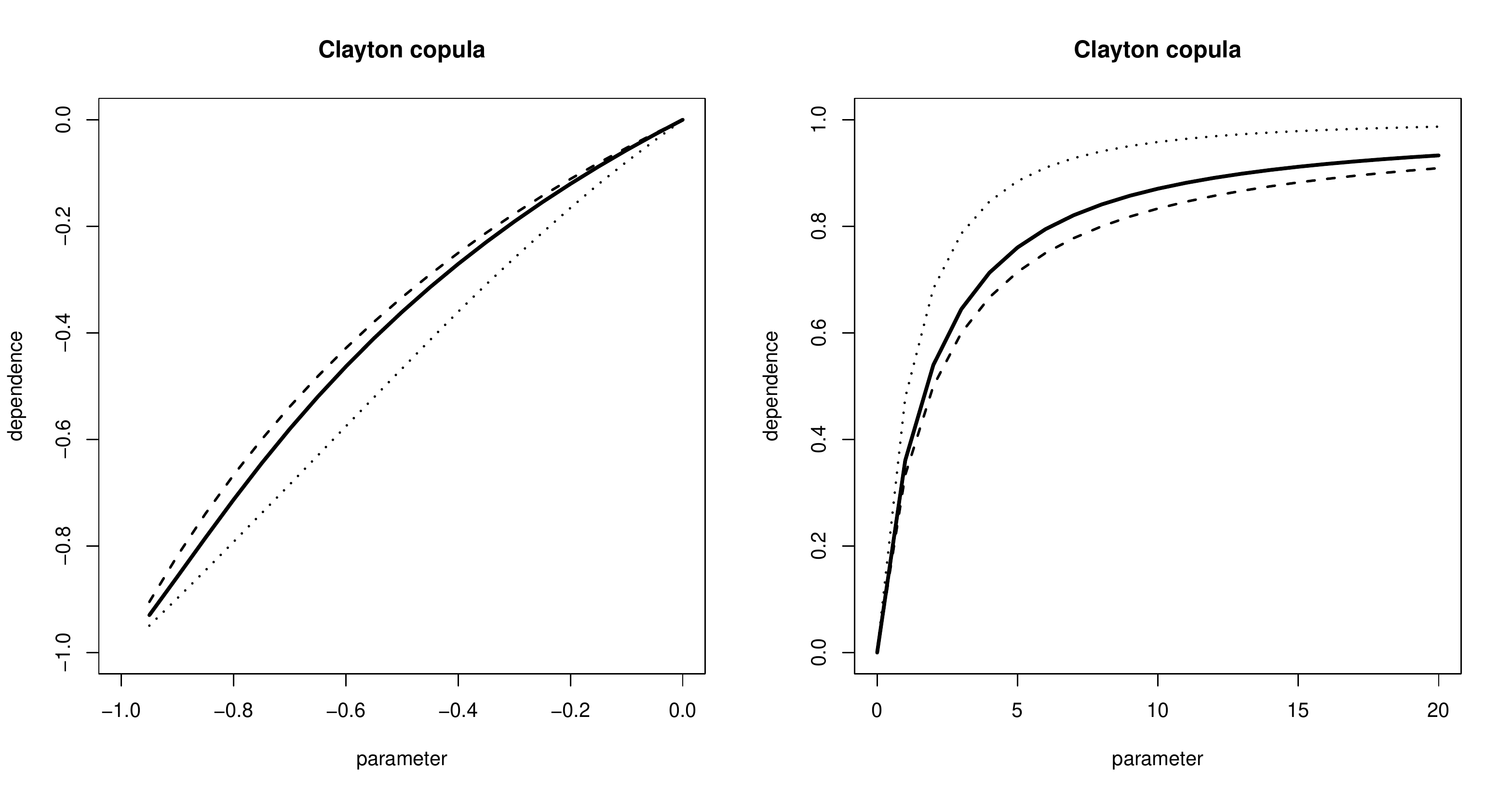}
\caption{Monotone dependence $\mu_{X,Y}$ (solid line), Kendall's $\tau_{X,Y}$ (dashed line), and Spearman's $\rho_{X,Y}$ (dotted line) in Example 3.}
\label{Ejemplo3}
\end{figure}
\end{center}

\noindent\textbf{Example 3} Let $(X,Y)$ be a random vector of continuous random variables. In this case the unique underlying subcopula in (\ref{Sklar2}) is, in fact, a copula, and therefore the value $\mu_{X,Y}$ will only depend on it, independently of the marginal distributions of $X$ and $Y.$ For continuous random variables, copula based concordance measures such as Kendall's $\tau_{X,Y}$ and Spearman's $\rho_{X,Y}$ are uniquely determined by:
\begin{equation*}
  \tau_{X,Y} = 4\!\int\!\!\!\int_{\mathbb{I}^{\,2}}\!C_{\theta}(u,v)\,dC_{\theta}(u,v)-1\,,\qquad\rho_{X,Y} = 12\!\int\!\!\!\int_{\mathbb{I}^{\,2}}\!C_{\theta}(u,v)\,dudv - 3\,,
\end{equation*}
and will be compared to the proposed monotone dependence $\mu_{X,Y}$ under the Clayton family of copulas:
\begin{equation*}
  C_{\theta}(u,v)\,=\,\big[\,\max(u^{-\theta}+v^{-\theta}-1, 0)\,\big]^{-1/\theta}\,,\qquad\theta\in[-1,+\infty[\,\setminus\{0\}
\end{equation*}
where $C_{-1}=W,$ $C_0=\Pi$ and $C_{\infty}=M.$ For Clayton copula with parameter $\theta$ it is possible to obtain explicitly $\tau_{X,Y}=\theta/(\theta+2),$ but for $\rho_{X,Y}$ a numerical approximation is required, which in this example is done by the \texttt{copula} R package by Hofert \textit{et al.}(2016). There is no explicit expression for $\mu_{X,Y}$ in this case, so to calculate it as in Corollary 2 a numerical maximization is performed applying the \texttt{nlm} function by the R Core Team (2016). See Figure \ref{Ejemplo3}.

\noindent The Clayton family of copulas was chosen as an example because it is \textit{comprehensive} (includes $W,$ $\Pi,$ and $M$) and therefore all values in the $[-1,1]$ interval may be reached for concordance measures with appropriate values of its parameter $\theta.$ As illustrated in Figure \ref{Ejemplo3} the behavior of $\mu_{X,Y}$ is similar to $\tau_{X,Y}$ and $\rho_{X,Y}.$

\section{Empirical subcopula}

Consider a bivariate random vector $(X,Y)$ where the random variables $X$ and $Y$ may be discrete, continuous or mixed type, not necessarily both of the same kind. Let $\{(x_1,y_1),\ldots,(x_n,y_n)\}$ denote a size $n$ sample of observations from $(X,Y).$ Since there may be repeated values among $\{x_1,\ldots,x_n\}$ let $\{r_1,\ldots,r_{m_1}\}$ be the set of distinct observed values of $X$ in the sample such that $r_1<\cdots<r_{m_1},$ where $m_1\leq n,$ and analogously let $\{s_1,\ldots,s_{m_2}\}$ be the ordered set of distinct observed values among $\{y_1,\ldots,y_n\},$ where $m_2\leq n.$

\bigskip

\noindent Let $p_{1i}$ be the proportion of the observed values of $X$ that are equal to $r_i$ and $p_{2j}$ the proportion of the observed values of $Y$ that are equal to $s_j,$ that is:
\begin{equation}\label{pli}
 p_{1i}\,:=\,\frac{1}{n}\sum_{k\,=\,1}^n \mathbf{1}\{x_k=r_i\}\,,\quad i\in\{1,\ldots, m_1\}\,, \qquad p_{2j}\,:=\,\frac{1}{n}\sum_{k\,=\,1}^n \mathbf{1}\{y_k=s_j\}\,,\quad j\in\{1,\ldots, m_2\}\,,
\end{equation}
where clearly $p_{1i}>0$ and $p_{2j}>0,$ and also $\sum_{i=1}^{m_1}p_{1i}=1$ and $\sum_{j=1}^{m_2}p_{2j}=1.$ Now define the sets $D_1=\{q_{10},q_{11},\ldots,q_{1m_1}\}$ and $D_2=\{q_{20},q_{21},\ldots,q_{2m_2}\}$ where $q_{10}=0=q_{20}$ and
\begin{equation}\label{qli}
  q_{1i}\,:=\,\sum_{k\,=\,1}^{i} p_{1k}\,,\quad i\in\{1,\ldots, m_1\}\,, \qquad q_{2j}\,:=\,\sum_{k\,=\,1}^{j} p_{2k}\,, \quad j\in\{1,\ldots, m_2\}\,,
\end{equation}
where clearly $0=q_{10}<q_{11}<\cdots<q_{1,m_1-1}<q_{1m_1}=1$ and $0=q_{20}<q_{21}<\cdots<q_{2,m_2-1}<q_{2m_2}=1.$ Then the set $D_1\times D_2$ is suitable as domain for a subcopula as in Definition 1. Let the function $S_n:D_1\times D_2\rightarrow\mathbb{I}$ be defined as $S_n(q_{10},q_{2j}):=0=:S(q_{1i},q_{20})$ for all $i$ and $j,$ and
\begin{equation}\label{subcopem}
  S_n(q_{1i},q_{2j})\,:=\,\frac{1}{n}\sum_{k\,=\,1}^n \mathbf{1}\{x_k\leq r_i\,,\,y_k\leq s_j\}\,,\quad i\in\{1,\ldots,m_1\}\,,\quad j\in\{1,\ldots,m_2\}\,,
\end{equation}
then it is straightforward to verify that $S_n$ is a subcopula and therefore we will call it \textit{empirical subcopula} associated to the observed sample $\{(x_1,y_1),\ldots,(x_n,y_n)\}.$ It should be noticed that the usual empirical joint distribution $F_n(r_i,s_j)=\frac{1}{n}\sum_{k=1}^n \mathbf{1}\{x_k\leq r_i\,,\,y_k\leq s_j\}=S_n(q_{1i},q_{2j}),$ that is $F_n$ and $S_n$ have the same range, but different domain since $F_n:\mathbb{R}^2\rightarrow\mathbb{I}.$ It is possible then to calculate (\ref{dS}) as:
\begin{equation}\label{dSn}
  d(S_n)\,=\,\max\{S_n(q_{1i},q_{2j})-q_{1i}q_{2j}:q_{1i}\in D_1,q_{2j}\in D_2\}\,-\,\max\{q_{1i}q_{2j}-S_n(q_{1i},q_{2j}):q_{1i}\in D_1,q_{2j}\in D_2\}
\end{equation}
and therefore a sample version of the monotone dependence measure would be $\mu(S_n)$ which is calculated accordingly to Definition 3 and (\ref{dSn}).

\bigskip

\noindent In the particular case where $X$ and $Y$ are both continuous random variables, an observed bivariate sample $\{(x_1,y_1),\ldots,(x_n,y_n)\}$ will contain non repeated values, and then $m_1=n=m_2,$ $r_i=x_{(i)}$ and $s_j=y_{(j)}$ (where $x_{(i)}$ stands for the $i$-th order statistic), $p_{1i}=\frac{1}{n}=p_{2j}$ for all $i$ and $j,$ $q_{1i}=\frac{i}{n}$ and $q_{2j}=\frac{j}{n},$ and $D_1=\mathbb{I}_{\,n}=D_2$ where $\mathbb{I}_{\,n}=\{0,\frac{1}{n},\ldots,\frac{n-1}{n}, 1\}.$ In this case the empirical subcopula $S_n:\mathbb{I}_{\,n}^{\,2}\rightarrow\mathbb{I}$ defined in (\ref{subcopem}) would be equivalent to:
\begin{equation}\label{SCcont}
  S_n\Big(\frac{i}{n}\,,\,\frac{j}{n}\Big)\,=\,\frac{1}{n}\sum_{k\,=\,1}^n\mathbf{1}\{x_k\leq x_{(i)}\,,\,y_k\leq y_{(j)}\}
\end{equation}
and $S_n(\frac{i}{n},0)=0=S_n(0,\frac{j}{n}),$ which agrees with the usual definition of \textit{empirical copula} for continuous random variables, see Definition 5.6.1 in Nelsen (2006). The expression ``empirical copula'' is somehow misleading since it is a subcopula with finite support $\mathbb{I}_{\,n}^{\,2},$ but not a copula. Of course, empirical subcopula (\ref{SCcont}) may be extended in a non-unique way to a copula, for example by bilinear interpolation as in Lemma 2.3.5 in Nelsen (2006), which is also known as \textit{checkerboard copula}, see Li \textit{et al.}(1997) or Durante and Sempi (2016).

\bigskip

\noindent Consequently, for observations from a pair of continuous random variables (\ref{dSn}) becomes:
\begin{equation}
  d(S_n)\,=\,\max\bigg\{S_n\Big(\frac{i}{n}\,,\,\frac{j}{n}\Big)-\frac{ij}{n^2}:i,j\in\{0,\ldots,n\}\bigg\}\,-\,\max\bigg\{\frac{ij}{n^2}-S_n\Big(\frac{i}{n}\,,\,\frac{j}{n}\Big):i,j\in\{0,\ldots,n\}\bigg\}
\end{equation}
and by elementary calculations:
\begin{equation}\label{cotasMonDep}
  -d(W_{S_n})\,=\,d(M_{S_n})\,=\,\left\{\begin{array}{cc}
	                                         \frac{1}{4} & \text{ if } n \text{ even,} \\
																						    { }    &    { }   \\
																				   \frac{n^2 - 1}{4n^2} & \text{ if } n \text{ odd, } \\
	                                       \end{array}\right.
\end{equation}
therefore by Definition 3:
\begin{equation}\label{empiricalMonDep}
  \mu(S_n)\,=\,\left\{\begin{array}{cc}
	                         4\hspace{0.3mm}d(S_n)(-1)^{\mathbf{1}\{d(S_n)\leq 0\}} & \text{ if } n \text{ even,} \\
																		      {  }   &   { }  \\
													 \frac{4n^2}{n^2 - 1}\hspace{0.3mm}d(S_n)(-1)^{\mathbf{1}\{d(S_n)\leq 0\}} & \text{ if } n \text{ odd, } \\
	                                                           \end{array}\right.
\end{equation} 

\medskip

\noindent where $\mathbf{1}\{d(S_n)\leq 0\}$ is equal to $1$ if $d(S_n)\leq 0$ and $0$ otherwise.

\bigskip

\noindent An R package \texttt{subcopem2D} by Erdely (2017) has been developed to perform the above calculations. \texttt{subcopem} function is for calculation of bivariate empirical subcopula matrix (\ref{subcopem}), induced partitions $D_1$ and $D_2,$ and $\mu(S_n)$ for a given bivariate sample of a pair of arbitrary type random variables. \texttt{subcopemc} function performs the same but it is specifically for a pair of continuous random variables with the possibility of faster calculations. \texttt{dependence} function calculates a matrix of pairwise dependence values for several variables. Examples are provided within the package.

\section{Conclusion}

A \textit{monotone dependence measure} $\mu_{X,Y}$ (Definition 3) is proposed for arbitrary type random variables $X$ and $Y$ based on the unique underlying subcopula given by \textit{Sklar's Theorem} (\ref{Sklar}), and its main properties are summarized in Theorem 1. Examples for discrete-discrete, continuous-discrete, and continuous-continuous pairs of random variables were analyzed, and in the particular case where both random variables are continuous $\mu_{X,Y}$ turns out to be a \textit{concordance measure} (Corollary 2), with the advantage that its definition is still unique in the general context of subcopulas, in contrast with other concordance measures that depend on non-unique extensions of copulas to subcopulas. Also, a sample version of the proposed dependence measure based on the \textit{empirical subcopula} has been provided, for the general case (\ref{dSn}) and for the particular case where both random variables are continuous (\ref{empiricalMonDep}), along with an R package \texttt{subcopem2D} by Erdely (2017) to perform such calculations.

\section*{Acknowledgement}

The author is specially grateful to the anonymous referees whose observations helped to significantly improve the present work. The present work was partially supported by Programa UNAM--DGAPA--PAPIIT \textbf{IN115817}.

\section*{References}

\noindent Cifarelli, D.M., Conti, P.L., and Regazzini, E. (1996) On the asymptotic distribution of a general measure of monotone dependence. \textit{Ann. Statist.} \textbf{24}, 1386--1399.
\medskip

\noindent Denuit, M. and Lambert, P. (2005) Constraints on concordance measures in bivariate discrete data. \textit{J. Mult. Anal.} \textbf{93}, 40--57.

\noindent Durante, F. and Sempi, C. (2016) \textit{Principles of Copula Theory.} CRC Press (Boca Raton).
\medskip

\noindent Erdely, A. (2017) \textit{subcopem2D: Bivariate Empirical Copula.} R package version 1.2 URL https://CRAN.R-project.org/package=subcopem2D
\medskip

\noindent Genest, C. and Ne\v{s}lehov\'a, J. (2007) A Primer on Copulas for Count Data. \textit{Astin Bull.} \textbf{37}(2), 475--515.
\medskip

\noindent Genest, C., Ne\v{s}lehov\'a, J., and R\'emillard, B. (2014) On the empirical multilinear copula process for count data. \textit{Bernoulli} \textbf{20}, 1344--1371.
\medskip

\noindent Hofert, M., Kojadinovic, I., Maechler, M., and Yan, J. (2016) \textit{copula: Multivariate Dependence with Copulas.} R package version 0.999-16 URL https://CRAN.R-project.org/package=copula
\medskip

\noindent Lehmann, E.L. (1966) Some concepts of dependence. \textit{Ann. Math. Statist.} \textbf{37}, 1137--1153.
\medskip

\noindent Li, X., Mikusi\'nski, P., Sherwood, H., and Taylor, M.D. (1997) On approximation of copulas. In Bene\v{s}, V. and \v{S}t\v{e}p\'an, J., editors, \textit{Distributions with given marginals and moment problems,} pp. 107--116. Kluwer (Dordrecht).
\medskip

\noindent Nelsen, R.B. (2006) \textit{An Introduction to Copulas.} Springer (New York).
\medskip

\noindent Ne\v{s}lehov\'a, J. (2007) On rank correlation measures for non-continuous random variables. \textit{J. Multivariate Analysis} \textbf{98}, 544--567.
\medskip

\noindent R Core Team (2016) \textit{R: A language and environment for statistical computing.} R Foundation for Statistical Computing (Vienna) URL https://www.R-project.org/
\medskip

\noindent Scarsini, M. (1984) On measures of concordance. \textit{Stochastica} \textbf{8} (3), 201--218.
\medskip

\noindent Schweizer, B. and Eolff, E.F. (1981) On nonparametric measures of dependence for random variables. \textit{Ann. Statist.} \textbf{9}, 879--885.
\medskip

\noindent Sklar, A. (1959) Fonctions de r\'epartition \`a $n$ dimensions et leurs marges. \textit{Publ. Inst. Statist. Univ. Paris,} \textbf{8}, 229--231.
\medskip

\end{document}